\nonstopmode
\documentclass[a4paper, 10pt]{amsart}
\usepackage{latexsym}
\usepackage{fancyhdr}
\usepackage{amsmath, amssymb}
\usepackage[ansinew]{inputenc}
\theoremstyle{plain}
\newtheorem*{definition*}{Definition}

\newtheorem*{lemma*}{Lemma}

\newtheorem*{theorem*}{Theorem}
\newtheorem{theorem}[subsection]{Theorem}
\newtheorem*{proposition*}{Proposition}
\newtheorem{proposition}[subsection]{Proposition}
\newtheorem*{corollary*}{Corollary}

\theoremstyle{remark}
\newtheorem*{remark*}{Remark}
\newtheorem{remark}[subsection]{Remark}
\def\cit#1#2{\ifx#1!\cite{#2}\else#2\fi} 
\def\fig#1#2{\ifx#1!\ref{fig:#2}\else\label{fig:#2}\fi} 
\def\ign#1{}             
\def\o{\circ}
\def\al{\alpha}
\def\be{\beta}
\def\ga{\gamma}
\def\de{\delta}

\def\vartheta{\theta}

\def\om{\omega}
\def\Ga{\Gamma}

\def\Om{\Omega}

\def\x{\times}
\def\g{\mathfrak g}
\def\p{\partial}
\def\on{\operatorname}

\def\today{\ifcase\month\or
 January\or February\or March\or April\or May\or June\or
 July\or August\or September\or October\or November\or December\fi
 \space\number\day, \number\year}
\title[A 2-cocycle on a group  of symplectomorphisms]
{A 2-cocycle on a group of symplectomorphisms}

\author[Ismagilov, Losik, Michor]{Rais S. Ismagilov, Mark Losik, Peter W. Michor}
\date{\today}
\thanks{ ML and PWM were supported by FWF Projekt P~17108-N04. 
PWM was supported by Centre Bernoulli, Lausanne}
\address{R.\ S.\ Ismagilov: Bauman Moscow State University, 2-nd Baumanskaya Str. 5, 
107005 Moscow, Russia.} 
\email{ismagil@serv.bmstu.ru}  
\address{M.\ Losik: Saratov State University, Astrakhanskaya 83,
410026 Saratov, Russia.}
\email{LosikMV@info.sgu.ru}
\address{P.\  W.\  Michor: Fakult\"at f\"ur Mathematik, Universit\"at Wien,
Nordbergstrasse 15, A-1090 Wien, Austria; {\it and}:
Erwin Schr\"odinger Institute of Mathematical Physics, Boltzmanngasse
9, A-1090 Wien, Austria.}
\email{Peter.Michor@esi.ac.at}

\begin{document}
\begin{abstract} 
For a symplectic manifold $(M,\om)$ with exact symplectic form 
we construct a 2-cocycle on the group of symplectomorphisms 
and indicate cases when this cocycle is not trivial. 
\end{abstract}
\maketitle

\section{Introduction} 
For a symplectic manifold $(M,\om)$  such that $H^1(M,\mathbb R)=0$ 
and the symplectic form $\om$ is exact we indicate 
a formula defining a 2-cocycle on the group $\on{Diff}(M,\om)$ of
symplectomorphisms with values in the trivial $\on{Diff}(M,\om)$-module
$\mathbb R$.
Let $G$ be a connected 
real simple Lie group and $K$ a maximal compact subgroup.
For the symmetric Hermitian space $M=G/K$ endowed with the induced 
symplectic structure, we prove that the restriction of this cocycle to the
group $G$ is non-trivial.
Thus this cocycle is non-trival on the whole group $\on{Diff}(M,\om)$, too. 
In particular, this implies that the cocycle is non-trivial for the
symplectic manifold $(\mathbb R^2\x M,\om_0+\om_M)$, where $(M,\om_M)$ is a
non-compact symplectic manifold with exact symplectic form $\om_M$ such
that $H^1(M,\mathbb R)=0$ and $\om_0$ is the standard symplectic form on
$\mathbb R^2$.

For the convenience of the reader, in an appendix we consider the
corresponding 2-cocycle on the Lie algebra of locally Hamiltonian and
Hamiltonian vector fields and indicate when this cocycle is non-trivial. 

Note that in \cite{Ism} a similar 2-cocycle was constructed for the group
of volume preserving diffeomorphisms on a compact $n$-dimensional manifold
$M$. This cocycle takes its values in the space $H^{n-2}(M,\mathbb R)$.
Neretin in \cite{Ne} constructed a 2-cocycle on the group of
symplectomorphisms with compact supports.

Throughout the paper $M$ is a connected $C^\infty$-manifold.

We thank Stefan Haller and Yurii Neretin for comments.

\section{Preliminaries}
We recall some standard facts on central extensions of groups and
two-dimen\-sio\-nal cohomology of groups (see, for example, \cite{ML}, ch. 4).

Consider a group $G$ and the field $\mathbb R$ as a trivial $G$-module. Let
$C^p(G,\mathbb R)$ be the set of maps from $G^p$ to $\mathbb R$ for $p>0$ and let
$C^0(G,\mathbb R)=\mathbb R$. Define a map  
$D^p:C^p(G,\mathbb R)\to C^{p+1}(G,\mathbb R)$ as follows:
for $f\in C^p(G,\mathbb R)$ and $g_1,\dots,g_{p+1}\in G$
\begin{multline}\label{D^p}
(D^pf)(g_1,\dots,g_{p+1})=f(g_2,\dots,g_{p+1})\\+
\sum_{i=1}^p(-1)^if(g_1,\dots,g_{i}g_{i+1},\dots,g_{p+1})+
(-1)^{p+1}f(g_1,\dots,g_p).
\end{multline}
By definition $C^*(G,\mathbb R)=(C^p(G,\mathbb R),D^p)_{p\ge 0}$ is the
standard complex of
nonhomogeneous cochains of the group $G$ with values in the
$G$-module $\mathbb R$ and its
cohomology $H^*(G,\mathbb R)=(H^p(G,\mathbb R))_{p\ge 0}$ is the cohomology
of the group $G$ with values in the trivial $G$-module $\mathbb R$.
Recall that a cochain $f\in C(G,\mathbb R)$ is called normalized
if $f(g_1,\dots,g_p)=0$ whenever at least one of the $g_1,\dots,g_p\in G$
equals the identity $e$ of $G$.
It is known that the inclusion of the subcomplex of normalized cochains
into $C^*(G,\mathbb R)$ induces an isomorphism in cohomology.

Let $f$ be a normalized $2$-cocycle of $G$ with values
in a trivial $G$-module $\mathbb R$.
Let $E(G,\mathbb R)=G\x\mathbb R$, with multiplication
$(g_1,a_1)(g_2,a_2)=(g_1g_2, a_1g_2+a_2+f(g_1,g_2))$ for $a_1,a_2\in
\mathbb R$
and $g_1,g_2\in G$.
Then $E(G,\mathbb R)$ is a group, and the natural
projection
$E(G,\mathbb R)=G\times\mathbb R\to G$ is a central extension of the group
$G$ by $\mathbb R$. The extension
$E(G,\mathbb R)$ is non-split iff the cocycle $f$ is non-trivial.

If $G$ is a topological group (finite-dimensional or infinite-dimensional
Lie group)
one can define a subcomplex $C^*_{\on{cont}}(G,\mathbb R)$
($C^*_{\on{diff}}(G,\mathbb R)$) of the complex $C^*(G,\mathbb R)$
(see \cite{Gui}, ch. 3) consisting of cochains which are continuous
(smooth) functions.
The cohomologies of the complexes
$C^*_{\on{cont}}(G,\mathbb R)$ and $C^*_{\on{diff}}(G,\mathbb R)$
are isomorphic whenever $G$ is a finite-dimensional Lie group
(see \cite{Gui}, ch. 3 and \cite{Mo}).
Note that if the 2-cocycle $f$ is continuous (differentiable), the
extension $E(G,\mathbb R)$ is a topological group (Lie group).

\section{A 2-cocycle on the group of symplectomorphisms}\label{2-cocycle}

Let $(M,\om)$ be a non-compact symplectic manifold such that $H^1(M,\mathbb R)=0$ 
and the symplectic form $\om$ is exact. 
Let $\om_1$ be a 1-form on $M$ such that $d\om_1=\om$.
Denote by $\on{Diff}(M,\om)$ the group of
symplectomorphisms of $M$.  We define a 2-cocycle on
the group $G=\on{Diff}(M,\om)$ with values in the trivial $G$-module
$\mathbb R$ as follows. Fix a point $x_0\in M$. Then for $g_1,g_2\in G$ we
put
\begin{equation}\label{C}
C_{x_0}(g_1,g_2)=\int_{x_0}^{g_2x_0}(g_1^*\om_1-\om_1),
\end{equation}
where the integral is taken along a smooth curve connecting the point $x_0$
with the point $g_2x_0$. Since $H^1(M,\mathbb R)=0$ the 1-form
$g_1^*\om_1-\om_1$ is exact and the value of this integral does not depend
on the choice of such a curve.

\begin{theorem}\label{cocycle} The function $C_{x_0}:G^2\to\mathbb R$
defined by (\ref{C}) is a normalized 2-cocycle on the group $G$ with values
in the trivial $G$-module $\mathbb R$. The cohomology class of $C_{x_0}$ is
independent of the choice of the point $x_0$ and the form $\om_1$.
\end{theorem}

\begin{proof}
By (\ref{D^p}) it is easy to check that $D^2C_{x_0}=0$. Moreover, the
2-cocycle $C_{x_0}$ is normalized.
Since for each $g\in G$ the 1-form $g^*\om_1-\om_1$ is exact,
for any points $x_1,x_2\in M$ we have $C_{x_1}-C_{x_2}=Da$, where $a$ is a
1-cochain on $G$ defined by $a(g)=\int_{x_1}^{x_2}(g^*\om_1-\om_1)$.
\end{proof}
By definition, the cocycle $C_{x_0}$ is a continuous function on $G\x G$. 

\begin{remark} Let $M$ be a manifold such that $H^1(M,\mathbb R)=0$ and 
let $\om$ be an exact 
2-form on $M$. Let $\on{Diff}(M,\om)$ be the group of diffeomorphisms of $M$ 
preserving the form $\om$. Then the formula (\ref{C}) for $g_1,g_2\in \on{Diff}(M,\om)$ 
gives a 2-cocycle on the group $\on{Diff}(M,\om)$  
and all statements of theorem (\ref{cocycle}) are true for this cocycle.
\end{remark}

Denote by $E(\on{Diff}(M,\om))$ the central extension of the group
$\on{Diff}(M,\om)$ by $\mathbb R$ defined by
the cocycle $C_{x_0}$.
Now we give a geometric interpretation of the extension
$E(\on{Diff}(M,\om))$.
We choose a form $\om_1$ with $d\om_1=\om$ 
and put $\om_2(g)=\int_{x_0}^x(\om_1-g^*\om_1)$.
Consider the trivial $\mathbb R$-bundle $M\x\mathbb R$.
Clearly, the
form $dt+\om_1$ is a connection with curvature $\om$ of this bundle.
Denote by $\on{Aut}(M\x\mathbb R,\om)$
the group of those bundle automorphisms which respect the
connection
$dt+\om_1$ and which are projectable to diffeomorphisms
in $\on{Diff}(M,\om)$. It is easy to check that the group
$\on{Aut}(M\x\mathbb R,\om)$
is isomorphic to the group
$E(\on{Diff}(M,\om))=\on{Diff}(M,\om)\times\mathbb R$ which acts as follows 
on $M\times \mathbb R$:
$(x,t)\to (g(x),\om_2(g)(x)+t+a)$, where $(x,t)\in M\x\mathbb R$ and
$(g,a)\in G\x\mathbb R$.
This gives an equivalent definition of the extension $E(\on{Diff}(M,\om))$
as a
group of automorphisms of the trivial principal
$\mathbb R$-bundle $M\times\mathbb R$ with connection $dt+\om_1$.

If we replace the form $\om_1$ by the form $\om_1+df$,
where $f$ is a smooth function on $M$, we get an
action of $G$ on $M\x\mathbb R$ which is related to the initial one
by the gauge transformation
$(x,t)\to (x,t-f(x))$ of the bundle $M\x\mathbb R\to M$.

\section{Examples of non-trivial $2$-cocycles}
\label{examples}

The authors are not able to prove that the cocycle $C_{x_0}$ is non-trivial 
for any symplectic manifold $M$ with an exact symplectic 2-form $\om$. 
In this section we prove that for some symplectic manifolds the
restrictions of this cocycle to some subgroups of
$G\subset\on{Diff}(M,\om)$ turn out to be non-trivial.

\subsection{The linear symplectic space $\mathbb R^{2n}$ and the Heisenberg
group}
\label{Heisenberg}

Consider the space $\mathbb R^{2n}$ with the standard symplectic form
$\om_0=\sum_{k=1}^n dx_k\wedge dx_{k+n}$ and the group $G=\mathbb R^{2n}$
acting on the space $\mathbb R^{2n}$ by translations.  Applying (\ref{C}) to
the form $\om_0$, the 1-form
$\om_1=\frac12\sum_{k=1}^n(x_{n+k}dx_k-x_kdx_{n+k})$, and the point
$x_0=0\in\mathbb R^{2n}$ we get a 2-cocycle on the group $G$ given by
$$
C_0(x,y)=\frac12\sum_{k=1}^n(x_ky_{n+k}-y_kx_{n+k}),
$$
where $x=(x_1,\dots,x_{2n})$ and $y=(y_1,\dots,y_{2n})$. The central
extension of the
group $\mathbb R^{2n}$ by $\mathbb R$ defined by this cocycle is the
Heisenberg group. This extension is non-split since the Heisenberg group is
noncommutative and thus the cocycle $C_0(x,y)$ is non-trivial.

\subsection{Symmetric Hermitian spaces and the Guichardet-Wigner
cocycle}\label{G-W}

Consider a non-compact symmetric space $M=G/K$, where $G$ is 
a connected real simple Lie group and where $K$ is a maximal 
compact subgroup. Then $M$ is diffeomorphic to $\mathbb R^n$, 
where $n=\dim M$. We suppose that $M$ 
admits a $G$-invariant complex structure, i.e., $M$ is a symmetric Hermitian 
space. This condition is satisfied (up to finite covering) 
for the following groups: $\on{SU}(p,q)$ ($p,q\ge 1$), $\on{SO}_0(2,q)$ 
($q=1$ or $q\ge 3$), 
$\on{Sp}(n,\mathbb R)$ ($n\ge 1$), $\on{SO}^*(2n)$ ($n\ge 2$), 
and certain real  
forms of $\on{E}_6$ and $\on{E}_7$.

Consider the symplectic manifold $(M,\om)$, where the symplectic form $\om$ 
is defined by the Hermitian metric on $M$. It is known that on each of the  
Lie groups mentioned above, 
in the complex $C_{\on{diff}}(G,\mathbb R)$ there is a non-trivial  
Guichardet-Wigner 2-cocycle (see \cite{G-W} and \cite{Gui}). By \cite{D-G} 
this cocycle is given as follows, up to a nonzero factor:
\begin{equation}\label{C_G} 
(g_1,g_2)\mapsto \int_{(x_0,g_1x_0,g_1g_2x_0)}\om,
\end{equation} 
where $g_1,g_2\in G$, $x_0=K\in G/K$, and the integral is taken over 
the oriented geodesic cone with vertex $x_0$ and the segment of a geodesic 
from $g_1x_0$ to $g_1g_2x_0$ as base. 

We prove that the restriction of the cocycle $C_{x_0}$ to the group $G$ 
is cohomologous to the cocycle given by (\ref{C_G}).   

For the base point $x_0$ we define a 1-cochain $\ga_{x_0}$ on the group 
$G$ as follows:
$$
\ga_{x_0}(g)=\int_{x_0}^{gx_0}\om_1,
$$
where $g\in G$ and the integral is taken along the geodesic segment 
from $x_0$ to $gx_0$.
Consider $C_{x_0}$ on $G$ given by 
formula (\ref{C}), where we choose for the curve between 
the points $x_0$ and $g_2x_0$ a geodesic segment from $x_0$ to $g_2x_0$. 
It is easy to check that on the group $G$ the cocycle  
$C_{x_0}+D\ga_{x_0}$ equals the cocycle given by (\ref{C_G}). 
Thus the cocycle $C_{x_0}$ on the group $G$  is non-trivial in the complex 
$C^*_{\on{diff}}(G,\mathbb R)$.

In particular, for the group $G=\on{SL}(2,\mathbb R)=\on{SU}(1,1)$ the
symmetric space $M=G/K$ is
the hyperbolic plane $H^2$ and $\om$ is the area form on $H^2$.
Instead of the group $\on{SL}(2,\mathbb R)$ we will later consider 
the group $\on{PSL}(2,\mathbb R)$ which acts effectively on $H^2$. Since
$\on{SL}(2,\mathbb R)$ is a
two-sheet cover of $\on{PSL}(2,\mathbb R)$, the cohomologies of these
groups with values in $\mathbb R$
are the same. It is easy to check that the corresponding symplectic
manifold $(M,\om)$ is
isomorphic to the symplectic manifold $(\mathbb R^2,\om_0)$, where $\om_0$
is the standard symplectic
form on $\mathbb R^2$. Unfortunately, for the groups 
$G\ne\on{SU}(1,1)$ mentioned above we do not know whether
the symplectic manifolds $(M,\om)$ and $(\mathbb R^{2n},\om_0)$, where
$\dim M=2n$, are isomorphic or not.

The following proposition is known. We do not know a good reference for this; 
then we give a short proof communicated to us by Yu.A. Neretin. 

\begin{proposition}\label{semisimple}
For each symmetric Hermitian space $M=G/K$, where $G$ is a connected simple 
Lie group and $K$ is its maximal compact subgroup, the corresponding Guichardet-Wigner 
cocycle is non-trivial in the complex $C^*(G,\mathbb R)$. 
\end{proposition}

\begin{proof} Let $p:\tilde G\to G$ be the universal cover and let $a=C_{x_0}$ 
be the Guichardet-Wigner cocycle for the group $G$.
Consider the corresponding to $a$ 
2-cocycle $\tilde a$ on $\tilde G$ induced by $p$. By construction, 
the cocycle $\tilde a$ is trivial, i.e., there is
a smooth function $b$ defined on $\tilde G$ such that for any $g,h\in\tilde G$ 
we have $\tilde a(g,h)=b(h)-b(gh)+b(g)$.
 
Assume that the cocycle $a$ is trivial in the complex $C^*(G,\mathbb R)$,
i.e., there exists a function $f:G\to\mathbb R$ such that for $g,h\in\Ga$
we have $a(g,h)=f(h)-f(gh)+f(g)$.

Then the difference $b-f\o p$ is a homomorphism $\tilde G\to\mathbb R$.
This homomorphism
vanishes near the identity element of $\tilde G$ since the group $\tilde G$
is simple, and thus it vanishes
on the whole of $\tilde G$ since $\tilde G$ is connected. Then the
function $f$ is smooth
and the cocycle $a$ is trivial in the complex $C^*_{\on{diff}}(G,\mathbb R)$.
This contradiction proves our statement. 
\end{proof}

\section {Cases of nontriviality of the cocycle $C_{x_0}$ for groups of
symplectomorphisms}
\label{nontrivial}

Let $(M,\om_M)$ be a non-compact symplectic manifold such that
$H^1(M,\mathbb R)=0$ with an exact symplectic form $\om_M$.
 
By formula (\ref{C}), the form $\om_M$
defines a 2-cocycle $C_{x_0}$ for the group $\on{Diff}(M,\om_M)$ 
with values in the trivial $\on{Diff}(M,\om_M)$-module $\mathbb R$.
The aim of this section is to indicate cases when this cocycle is
non-trivial and thus the corresponding central extension of the group 
$\on{Diff}(M,\om_M)$ by $\mathbb R$ is non-split.

Let $M=G/K$ be an Hermitian symmetric space $M$ and let $(M,\om)$ be the corresponding 
symplectic manifold  which we considered in subsection \ref{G-W}.  

\begin{theorem}\label{G/K}
For the Hermitian symmetric space $M=G/K$ and for the corresponding
symplectic manifold $(M,\om)$
the cocycle $C_0$ on the group $\on{Diff}(M,\om)$ is non-trivial. 
\end{theorem}

\begin{proof} Since the group $G$ is a subgroup of the group
$\on{Diff}(M,\om)$ the statement follows
from proposition \ref{semisimple}.
\end{proof}

Recall that the symplectic manifold $(H^2,\om)$ where $\om$ is the area
form is symplectomorphic to $(\mathbb R^2,\om_0)$ where $\om_0$ is the
standard symplectic form.

\begin{theorem} \label{R^2xM}
Let $(M,\om)$ be a non-compact symplectic manifold such that the symplectic form $\om_M$ 
is exact and let $H^1(M,\mathbb R)=0$.
Consider the product $\mathbb R^2\times M$
of the manifold $\mathbb R^2$ and $M$ as a symplectic manifold with the symplectic form 
$\om=\om_0 + \om_M$.
Then for each point $x_0\in\mathbb R^2\times M$ the cocycle $C_{x_0}$ on the group 
$\on{Diff}(\mathbb R^2\x M,\om)$ is non-trivial.
\end{theorem}

\begin{proof}  
Choose $\om_{M,1}\in\Om^1(M)$ with
$d\om_{M,1}=\om_M$ and let $\om_1=x\,dy +\om_{M,1}$.
The group $\on{Diff}(\mathbb R^2,\om_0)$
acting on the first factor $\mathbb R^2$ of $\mathbb R^2\times M$
is naturally included as a subgroup into
the group $\on{Diff}(\mathbb R^2\times M,\om)$.
Thus $g^*\om_1-\om_1= g^*(x\,dy)-x\,dy$ for all $g$ in the subgroup
$\on{Diff}(\mathbb R^2,\om_0)$.
Thus the cocycle $C_{x_0}$ constructed from the
form $dx\wedge dy +\om_M$ on $\mathbb R^2\times M$ restricts to
a nontrial cocycle on the subgroup of $\on{Diff}(\mathbb R^2,\om_0)$
by proposition \ref{semisimple} applied to the group $\on{PSL}(2,\mathbb R)$.
\end{proof}
We leave to the reader to formulate the corresponding results for 
other symmetric Hermitian spaces $G/K$ instead of $H^2$.  

\subsection{Problem}
Consider an open disk $M$ in the Euclidean plane equipped with the
standard area 2-form $\om$. Is the 2-cocycle $C_{x_0}$ defined by the 
form $\om$ non-trivial?

\section{appendix}

In this appendix, for a symplectic manifold $(M,\om)$ we define a 2-cocycle
on the Lie algebra $\on{Vect}(M,\om)$ of locally Hamiltonian or Hamiltonian  
vector fields, corresponding to the 2-cocycle $C_{x_0}$ 
on the group $\on{Diff}(M,\om)$,  
and study conditions of its nontriviality.

Let $\g$ be a Lie algebra over $\mathbb R$ and let $\mathbb R$ be the
trivial $\g$-module. Denote by $C^p(\g,\mathbb R)$ the space of
skew-symmetric
$p$-forms on $\g$ with values in $\mathbb R$. 
For $c\in C^p(\g,\mathbb R)$ and $x_1,\dots,x_{p+1}$ put
\begin{equation}\label{de}
(\de^p c)(x_1,\dots,x_{p+1})=\sum_{i<j}(-1)^{i+j}c([x_i,x_j],x_1,\dots,\hat
x_i,\dots,\hat x_j\dots,x_{p+1}),
\end{equation}
where, as usual, $\hat x$ means that $x$ is omitted.
Then $C^*(\g,\mathbb R)=(C^p(\g,\mathbb R),\de^p)_{p\ge 0}$ is the complex
of standard cochains of the Lie algebra $\g$ with values in the trivial
$\g$-module $\mathbb R$ and
the cohomology $H^*(\g,\mathbb R)$ of this complex is the cohomology of the
Lie algebra $\g$ with values in the trivial $\g$-module $\mathbb R$.

In particular, there is a bijective correspondence between 
$H^2(\g,\mathbb R)$ and the set of isomorphism classes of central 
extensions of the Lie algebra $\g$ by $\mathbb R$.

Let $(M,\om)$ be a symplectic manifold. Denote by
$\on{Vect}(M,\om)$ the Lie algebra of locally Hamiltonian vector fields   
and by $\on{Vect}_0(M,\om)$ the Lie algebra of Hamiltonian vector fields on $M$. 
For a point $x_0\in M$ and $X,Y\in \on{Vect}(M,\om)$ put
$c_{x_0}(X,Y)=\om(X,Y)(x_0)$.

\begin{proposition}\label{c_0}
The function $c_{x_0}:\g^2\to\mathbb R$ is a 2-cocycle on the Lie algebra
$\g$ with values in the trivial $\g$-module $\mathbb R$. The cohomology
class of $c_{x_0}$ is independent of the choice of the point $x_0$. 
\end{proposition}

\begin{proof} The proof is given by direct calculations and is based on 
the standard formulas $[\mathcal L_X,\mathbf i_Y]=\mathbf i_{[X,Y]}$  and
$\mathcal L_X=\mathbf i_Xd+d\mathbf i_X$, where $\mathbf i_X$ is the operator
of the inner product by $X$ and  $\mathcal L_X$ is the Lie derivative with respect 
to a vector field $X$, (see, for example, \cite{Go}, ch. 4). In particular,
we have for any $x\in M$ and $X,Y\in\on{Vect}(M,\om)$ the following equality   
\begin{equation}\label{c_x}
c_x(X,Y)-c_{x_0}(X,Y)=-\int_{x_0}^x\mathbf i_{[X,Y]}\om.
\end{equation}
\end{proof}

Let $G$ be a Lie group and let $\g$ be its Lie algebra. We have a natural
homomorphism
of complexes $C^*_{\on{diff}}(G,\mathbb R)\to C^*(\g,\mathbb R)$ (see, for
example, \cite{Gui}, ch. 3). In particular, if $c\in
C^2_{\on{diff}}(G,\mathbb R)$, the corresponding cochain 
$\tilde c\in C^2(\g,\mathbb R)$ is defined as follows:
$$
\tilde c(X,Y)=\frac{\p^2}{\p t\p s}(c(\on{exp}tX,\on{exp}sY)-
c(\on{exp}sY,\on{exp}tX))_{t=0,s=0}
$$
where $X,Y\in\g$.

Let $G$ be a Lie group of diffeomorphisms of $M$ contained in the group
$\on{Diff}(M,\om)$. Then
for the 2-cocycle $c=C_{x_0}$ of section \ref{2-cocycle}, the cocycle
$\tilde c$ is cohomologous to the restriction of the cocycle $c_{x_0}$ to
the Lie algebra $\g$ of $G$. Unfortunately, we cannot apply this procedure
to the whole group
$\on{Diff}(M,\om)$ and the Lie algebra $\on{Vect}(M,\om)$. Therefore, the
problems of
nontriviality of 2-cocycles $C_{x_0}$ on the group $\on{Diff}(M,\om)$ and
$c_{x_0}$ on
the Lie algebra $\on{Vect}(M,\om)$ should be solved independently.

For each $X\in\on{Vect}(M,\om)$ denote by $\al_X$ the closed 1-form such that 
$\al_X=\mathbf i_X\om$.  
For all vector fields $X,Y\in\on{Vect}(M,\om)$ we have the following
equality:
\begin{equation}\label{om}
\om(X,Y)\om^n=n\al_X\wedge\al_Y\wedge\om^{n-1}
\end{equation}
which can be easily checked in Darboux coordinates.

Denote by $X_f$ a Hamiltonian vector field defined by a function $f\in
C^\infty(M)$.
Consider the Poisson algebra $\on{P}(M)=\on{P}(M,\om)$ on $(M,\om)$, i.e., the
algebra $C^\infty(M)$ endowed with the Poisson
bracket $\{f,g\}=-\om(X_f,X_g)$ for $f,g\in C^\infty(M)$. 
 
The map $\on{P}(M)\to\on{Vect}_0(M,\om)$ given by $f\to X_f$ is a
homomorphism of Lie algebras which defines an extension of
$\on{Vect}_0(M,\om)$ by $\mathbb R$. It is easy to check that this extension 
is isomorphic to one given by the cocycle $-c_{x_0}$.

\begin{theorem} For a non-compact symplectic manifold $(M,\om)$ the cocycle
$c_{x_0}$ on the Lie algebras $\on{Vect}(M,\om)$ and $\on{Vect}_0(M,\om)$
is non-trivial.
\end{theorem}

\begin{proof} It suffices to prove our statement for the Lie algebra 
$\on{Vect}_0(M,\om)$.

First we prove that for each form $\be\in\Om^{2n-1}(M)$ there is a unique
form $\al\in\Om^1(M)$ such that $\be=\al\wedge\om^{n-1}$. Indeed, using
Darboux coordinates it is easy to check that this has a unique local
solution $\al$.
These are compatible and we get a global solution by gluing them.

Note that for each form $\al\in\Om^1(M)$ there is a positive integer $N$
and $2N$ functions $f_k,g_k\in C^\infty(M)$ $(k=1,\dots,N)$ such that
$\al=\sum_{k=1}^Nf_kdg_k$ which follows easily from the existence (by
dimension theory) of a finite atlas for $M$.

Since $H^{2n}(M,\mathbb R)=0$ there is a form $\be\in\Om^{2n-1}(M)$ such
that
$\om^n=d\be$. Then we have 
$\om^n=\sum_{k=1}^Ndf_k\wedge dg_k\wedge\om^{n-1}$.
By (\ref{om}) and using this equality we get
\begin{equation}\label{-n}
\sum_{k=1}^N\{f_k,g_k\}=-n.
\end{equation}
Assume that the extension $P(M)\to\on{Vect}_0(M,\om)$
$P(M)\to\on{Vect}_0(M,\om)$ is split. Then $P(M)$ is a direct
sum of the space of constant functions on $M$ and an ideal isomorphic to
$\on{Vect}_0(M,\om)$
by $P(M)\to\on{Vect}_0(M,\om)$. Equality (\ref{-n}) means that these
summands have nonzero
intersection. This contradiction proves the statement.  
\end{proof}

Now we consider a compact symplectic manifold $(M,\om)$. It is known that 
the extension $P(M)\to\on{Vect}_0(M,\om)$ is split. 
 
For a closed form $\al$ denote by $[\al]$ the cohomology class of $\al$.
Denote by $L$ the linear map $H^p(M,\mathbb R)\to H^{p+2}(M,\mathbb R)$
defined
by $a\to a\smallsmile [\om]$, where $a\in H^p(M,\mathbb R)$.

\begin{theorem} Let $(M,\om)$ be a compact symplectic manifold.
The cocycle $c_{x_0}$ on the Lie algebra $\on{Vect}(M,\om)$ is non-trivial
iff the linear map
$$
L^{n-1}:H^1(M,\mathbb R)\to H^{2n-1}(M,\mathbb R)
$$ 
is not equal zero.
\end{theorem}
\begin{proof} We may assume that $\int_M\om^n=1$. Put for brevity
$V=\on{Vect}(M,\om)$ and $V_0=\on{Vect}_0(M,\om)$.  Set 
$$
b(X,Y)=\int_M\om(X,Y)\om^n,
$$
where $X,Y\in V$. It is easy to check that $b$ is a 2-cocycle on $V$.

Multiplying both sides of equality (\ref{c_x}) by $\om^n$ and integrating over
$M$ we get
\begin{equation}\label{b=c}
b(X,Y)-c_{x_0}(X,Y)=\int_M\left(\int_{x_0}^x\mathbf
i_{[X,Y]}\om\right)\om^n.
\end{equation}
Since the right hand side of (\ref{b=c}) is a coboundary of a 1-cochain
in $C^1(V,\mathbb R)$, the cocycles $c_{x_0}$ and $b$ are cohomologous.
By (\ref{om}) we have
\begin{equation}\label{b}
b(X,Y)=n\int_M\al_X\wedge\al_Y\wedge\om^{n-1},
\end{equation}
for any $X,Y\in V$.
If $X\in V_0$ the form $\al_X$ is exact, and $b(X,Y)=0$ by (\ref{b}). This
proves (1).

Suppose that the cocycle $b$ is trivial, i.e., there is a linear functional
$f$ on
$V$ such that for any $X,Y\in V$ we have $b(X,Y)=f([X,Y])$.
By \cite{Ar} we have $[V,V]=[V_0,V_0]=V_0$. This implies $b=0$.
So the cocycle $b$ is trivial iff it equals zero. By (\ref{b}) and the
Poincar\'e duality this implies that $L^{n-1}=0$ on $H^1(M,\mathbb R)$.
This proves (2).
\end{proof}

We know no example when $H^1(M,\mathbb R)\ne 0$ and the map $L^{n-1}=0$.
Moreover, if $M$ is a compact K\"ahlerian manifold the map
$L^{n-1}:H^1(M,\mathbb R)\to H^{2n-1}(M,\mathbb R)$ is an isomorphism (see,
for example,
\cite{We}, ch.\ 4). Thus in this case the cocycle $c_{x_0}$ is non-trivial
whenever $H^1(M,\mathbb R)\ne 0$.

\end{document}